\documentclass[10pt]{article}
\usepackage{amsmath,amssymb,graphicx}

\begin{document}

\title{Phase Lag Sensitivity Analysis for Numerical Integration}

\author{D. S. Vlachos\thanks{dvlachos@uop.gr}, Z. A. Anastassi\thanks{zackanas@uop.gr} and T. E. Simos
\thanks{Highly Cited Researcher, Active Member of the European Academy of Sciences and Arts, Address: Dr. T.E. Simos, 26 Menelaou Street, Amfithea - Paleon Faliron, GR-175 64 Athens, GREECE, Tel: 0030 210 94 20 091, e-mail: tsimos.conf@gmail.com, tsimos@mail.ariadne-t.gr}}

\maketitle

\begin{center}
Laboratory of Computational Sciences,\\
Department of Computer Science and Technology,\\
University of Peloponnese,\\
GR-22 100, Terma Karaiskaki, Tripolis, Greece\\
~
\date{\today}
\end{center}

\begin{abstract}
In the field of numerical integration, methods specially tuned on oscillating
functions, are of great practical importance. Such methods are needed in various
branches of natural sciences, particularly in physics, since a lot of physical
phenomena exhibit a pronounced oscillatory behavior. Among others, probably the
most important tool used to construct efficient methods for oscillatory problems
is the exponential (trigonometric) fitting. The basic characteristic of these methods
is that their phase lag vanishes at a predefined frequency. In this work, we
introduce a new tool which improves the behavior of exponentially fitted numerical methods. The new technique is based on the vanishing of the first derivatives of the phase lag function at the fitted frequency. It is proved in the text that these methods present improved characteristics in oscillatory problems.
\end{abstract}

\emph{PACS: }{0.260, 95.10.E}

\section{Introduction}
The numerical integration of systems of ordinary differential equations with
oscillatory solutions has been the subject of research during the past decades.
This type of ODEs is often met in real problems, like the Schr\"{o}dinger equation
and the N-body problem. For problems having highly oscillatory solutions, the standard non-specialized
methods can require a huge number of steps to track the oscillations. One way to obtain a more efficient integration process is to construct numerical
methods with an increased algebraic order, although the implementation of high algebraic order methods is not evident.

On the other hand, there are some special techniques for optimizing numerical methods. Trigonometrical fitting and phase-fitting are some of them, producing methods with
variable coefficients, which depend on $v = \omega h$, where $\omega$ is the dominant frequency of the problem and $h$ is the step length of integration. More precisely, the coefficients of a general linear method are found from
the requirement that it integrates exactly powers up to degree $p+1$. For problems having oscillatory
solutions, more efficient methods are obtained when they are exact for every linear
combination of functions from the reference set
\begin{equation}
\{1, x, \ldots , x^K , e^{\pm \mu x},\ldots , x^P e^{\pm \mu x}\}\label{equ_exp_fit}
\end{equation}
This technique is known as exponential (or trigonometric if $\mu=i\omega$) fitting and has a long history \cite{gautschi_NM_3_381_61}, \cite{lyche_NM_19_65_72}. The
set (\ref{equ_exp_fit}) is characterized by two integer parameters, $K$ and $P$ . The set in which there
is no classical component is identified by $K =-1$ while the set in which there is no
exponential fitting component (the classical case) is identified by $P =-1$. Parameter
$P$ will be called the level of tuning. An important property of exponential fitted algorithms is that they tend to the classical ones when the involved frequencies tend to zero, a fact which allows to say that exponential fitting represents a natural extension of the classical polynomial fitting. The examination of the convergence
of exponential fitted multistep methods is included in Lyche's theory \cite{lyche_NM_19_65_72}. There is a large number of significant methods presented with high practical importance thats have been presented in the bibliography (see for example \cite{simos_Book_CMATV1_RCS_00}, \cite{chawla_JCAM_15_329_86}, \cite{raptis_CPC_14_1_78}, \cite{anastassi_IJMPC_15_1_04}, \cite{anastassi_MCM_42_877_05}, \cite{anastassi_JMC_41_79_07}, \cite{anastassi_NA_10_301_05}, \cite{lambert_JIMA_18_189_76}, \cite{cash_JNAIAM_1_81_06}, \cite{iavernaro_JNAIAM_1_91_06}, \cite{mazzia_JNAIAM_1_131_06}, \cite{berghe_JNAIAM_1_241_06}, \cite{psihoyios_CL_2_51_06}, \cite{simos_CL_1_37_05}, \cite{simos_CL_1_45_07}. The general theory is presented in detail in \cite{ixaru_Book_EF_KAP_04}.

Considering the accuracy of a method when solving oscillatory problems,
it is more appropriate to work with the phase-lag, rather than the usually used principal local truncation error. We mention
the pioneering paper of Brusa and Nigro \cite{brusa_IJNME_15_685_80}, in which the phase-lag property was
introduced. This is actually another type of a truncation error, i.e. the angle between
the analytical solution and the numerical solution. On the other hand, exponential fitting is accurate only when a good estimate of the dominant frequency of the solution
is known in advance. This means that in practice, if a small change in the dominant frequency is introduced, the efficiency of the method can be dramatically altered. It is well known, that for equations similar to the harmonic oscillator, the most efficient exponentially fitted methods are
those with the highest tuning level. In the case of the Schr\"{o}dinger equation, this result
was already obtained for particular two- and four-step exponentially fitted multistep methods based on
an expensive error analysis, see for example \cite{ixaru_CPC_19_23_80},\cite{simos_JCAM_30_251_90}, \cite{simos_IMAJNA_11_347_91} and \cite{simos_PLA_177_345_93}.

In this paper we present a methodology for optimizing numerical methods,
through the use of phase-lag function and its derivatives with respect to $v$. More specifically, given a classical (that is with constant coefficients) numerical method, we
can provide a family of optimized methods, each of which has zero phase lag (the case of trigonometric fitting) or
zero $PL$ and $PL'$ or zero $PL$, $PL'$ and $PL''$ etc.
With this new technique we provide methods that perform well during the
integration of the Schr\"odinger equation for high values of energy, but also that
perform well on other real problems with oscillatory solution, like the N-body
problem.

\section{Phase lag analysis}
Below we will consider for simplicity only first order differential
equations, although the same results can be easily obtained for second order
equations too. Consider the test problem
\begin{equation}
\frac{dy(t)}{dt}=i\omega _0 y(t),\;y(0)=1
\label{equ_test_equ}
\end{equation}
with exact solution
\begin{equation}
y(t)=e^{i\omega _0 t}
\end{equation}
where $\omega _0$ is a non-negative real value. Let $\hat{\Phi}(h)$ be a numerical
map which when it is applied to a set of known past values, it produces a numerical estimation of
$y(t+h)$. If we assume that all past values are known exactly, then the
numerical estimation $\hat{y}(t+h)$ of $y(t+h)$ will be
\begin{equation}
\hat{y}(t+h)=\alpha (\omega _0 h) \cdot e^{i(\omega _0 t+\phi (\omega _0 h) )}
\end{equation}
while the exact solutiion is $e^{i (\omega _0 t +\omega _0 h)}$. Then
\begin{equation}
L=\frac{\hat{y}(t+h)}{y(t+h)}=\alpha (\omega _0 h)e^{-i(\omega _0 h-\phi (\omega
_0 h) ) }
\label{equ_def_PL}
\end{equation}
In the above equation (\ref{equ_def_PL}), the term $\alpha (\omega _0 h)$ is
called the \textit{amplification factor}, while the term $l(\omega _0 h)=\omega
_0 h-\phi(\omega _0 h)$ is called the \textit{phase lag} of the numerical map.
In the case that $\alpha (\omega _0 h)=1$ and $l(\omega _0 h)=0$, we say that
the numerical map $\Phi(h)$ is \textit{exponentially fitted} at the frequency
$\omega _0$ and at the step size $h$.
\par Suppose now that the method $\hat{\Phi}(h)$ has been designed in order to
solve exactly equation (\ref{equ_test_equ}). But in practice, only an estimation of the frequency $\omega _0$ is known. Thus, it is of great importance to know the behavior of the method at frequencies close to the estimated one, so we apply the method to the equation
\begin{equation}
\frac{dy(t)}{dt}=i\omega y(t),\;y(0)=1
\end{equation}
and calculate the phase lag $l(u)$, where $u=\omega h$. Since the method integrates exactly equation (\ref{equ_test_equ}), the phase lag function $l(u)$
has a zero at point $u_0=\omega _0 h$ and is given by
\begin{equation}
l(u)=u-\phi(u)
\end{equation}
But since $\phi (u_0)=u_0$ we have
\begin{equation}
l(u)=u-\phi (u_0)-\frac{d\phi(u)}{du}|_{u=u_0}(u-u_0)+\sum_{n=2}^{\infty}\frac{d^{(n)}\phi(u)}{du^n}|_{u=u_0}\frac{(u-u_0)^n}{n!}
\end{equation}
or
\begin{equation}
l(u)=(u-u_0)\left(1-\frac{d\phi(u)}{du}|_{u=u_0}\right)+\sum_{n=2}^{\infty}\frac{d^{(n)}\phi(u)}{du^n}|_{u=u_0}\frac{(u-u_0)^n}{n!}
\end{equation}
and thus we conclude that in order to maximize the phase lag order at
frequencies close to $u_0=\omega _0 h$ we must at least have
\begin{equation}
\frac{d\phi(u)}{du}|_{u=u_0}=1
\end{equation}
while in order to obtain even higher order in phase lag, the higher derivatives
of $\phi(u)$ at point $u_0$ must vanish.Since
\begin{equation}
l'(u)=1-\phi '(u),\;l^{(p)}(u)=\phi ^{(p)}(u),\;p>1
\end{equation}
we have the following theorem
\par\textbf{Theorem 1.} Consider a linear method $\hat{\Phi}$ which solves exactly the equation (\ref{equ_test_equ}) and when it is applied to the equation $y'=i\omega y$, with $\omega=\omega _0 +\delta$, it produces a phase lag function $l(u)$, $u=\omega h$. Then, if the phase lag function has its $s$ first derivatives at point $v_0=\omega _0 h$ equal to zero, then the phase lag function $l(u)$ is of order at least $s$ in $\delta$.
\par Supose now that the method $\hat{\Phi}$ depends on $M$ independent parameters and let $\hat{ \Phi _c}$ the classical method which is constructed by setting $M$ equations maximizing the algebraic order of the method. Let $\hat{\Phi _s}$ be the method which is constructed with $M-1-s$ equations maximizing the algebraic order, $1$ equation for vanishing the phase lag at a frequency $\omega _0 h$ and $s$ equations for vanishing the $s$ first derivatives of the phase lag at the same point $\omega _0 h$. It is easily now calculated that the local truncation error of the method $\hat{\Phi _s}$ is given by
\begin{equation}
lte_s=\alpha (u) \left(e^{-il(u)}-1\right)
\end{equation}
when the working frequency $\omega _0 \rightarrow 0$. But since the method is at least of order $M-1-s$ we have
\begin{equation}
\begin{array}{l}
lte_s=h^{M-s}(l(u_0)-il'(u_0)(\omega -\omega _0)h\\
+\sum_{j=2}^{\infty}\left(\lambda _j l^{(j)}(u_0)+g_j(l^{(1)}(u_0),\ldots,l^{(j-1)}(u_0))\right)h^j
\end{array}
\end{equation}
with $g_j$ be a polynomial function of the $j-1$ first derivatives of $l(u)$ at point $u_0$ with $g_j(0,0,\ldots,0)=0$. Since now $l(u_0)=0$ and the $s$ first derivatives of $l(u)$ at $u_0$ vanish, we have that in the limit $\omega _0\rightarrow 0$
\begin{equation}
lte_s=ch^{M+1}=lte_c
\end{equation}
where $lte_c$ is the local truncation error of the classical method $\hat{\Phi} _c$. Thus we have the following theorem
\par\textbf{Theorem 2.} Consider a linear method $\hat{\Phi}_s$ which is constructed by demanding (i) maximal algebraic order and (ii) that the phase lag and its $s$ first derivatives vanish at some given frequency $\omega _0$. Then, when $\omega _0 \rightarrow 0$, the method is identical with the classical one $\hat{\Phi} _c$ which is constructed by demanding only maximal algebraic order.

\section{Numerical results}
In order to follow the dynamics of the method constructed, vanishing the first derivatives of the phase lag function, we consider the simple $2-step$ symmetric formula
\begin{equation}
y_{n-1}-2\cdot y_n+y_{n+1}=h^2\left(b_0\cdot(f_{n-1}+f_{n+1})+b_1f_n\right)
\label{equ_method}
\end{equation}
for the solution of the $2nd-order$ equation
\begin{equation}
\frac{d^2y}{dt^2}=f(y)
\end{equation}
The coefficients of the method are calculated in three different cases as follows:
\begin{enumerate}
\item $b^c$ are the coefficients for the classical method, where only the maximization of the algebraic order is taken into account
\item $b^t$ are the coefficients for the method where trigonometric fitting in a frequency $\omega$ ($v=\omega h$) is taken into account and
\item $b^s$ are the coefficients for the method where both trigonometric fitting and vanishing of the first derivative of the phase lag function is taken into account.
\item $b^{sd}$ are the coefficients for the method where both trigonometric fitting and vanishing of the first and second derivatives of the phase lag function is taken into account. In order to obtain this, the coefficient of $y_n$ in equation (\ref{equ_method} is perturbed from $-2$ to $-2+a(v)$.
\end{enumerate}
The coefficients for the three cases are given
\begin{equation}
b^c_0=\frac{1}{12},\;b^c_1=\frac{5}{6}
\end{equation}
\begin{eqnarray}
b^t_0=&&{\frac {1}{12}}+{\frac {1}{120}}{v}^{2}+{\frac {17}{20160}}{v}^{4}+{
\frac {31}{362880}}{v}^{6}+{\frac {691}{79833600}}{v}^{8}+\nonumber \\&&{\frac {5461
}{6227020800}}{v}^{10}+{\frac {929569}{10461394944000}}{v}^{12}+O
 \left( {v}^{14} \right)
\end{eqnarray}
\begin{eqnarray}
b^t_1=&&{\frac {5}{6}}-{\frac {1}{60}}{v}^{2}-{\frac {17}{10080}}{v}^{4}-{
\frac {31}{181440}}{v}^{6}-{\frac {691}{39916800}}{v}^{8}-\nonumber \\&&{\frac {5461
}{3113510400}}{v}^{10}-{\frac {929569}{5230697472000}}{v}^{12}+O
 \left( {v}^{14} \right)
\end{eqnarray}
\begin{eqnarray}
b^s_0=&&{\frac {1}{12}}+{\frac {1}{120}}{v}^{2}+{\frac {17}{20160}}{v}^{4}+{
\frac {31}{362880}}{v}^{6}+{\frac {691}{79833600}}{v}^{8}+\nonumber \\&&{\frac {5461
}{6227020800}}{v}^{10}+{\frac {929569}{10461394944000}}{v}^{12}+O
 \left( {v}^{14} \right)
\end{eqnarray}
\begin{eqnarray}
b^s_1=&&{\frac {5}{6}}-{\frac {1}{60}}{v}^{2}+{\frac {5}{2016}}{v}^{4}+{
\frac {29}{181440}}{v}^{6}+{\frac {139}{7983360}}{v}^{8}+\nonumber \\&&{\frac {5459}
{3113510400}}{v}^{10}+{\frac {185917}{1046139494400}}{v}^{12}+O
 \left( {v}^{14} \right)
\end{eqnarray}
\begin{eqnarray}
b^{sd}_0=&&{\frac {1}{12}}+{\frac {1}{80}}{v}^{2}+{\frac {41}{20160}}{v}^{4}+{
\frac {1219}{3628800}}{v}^{6}+{\frac {8887}{159667200}}{v}^{8}+\nonumber \\&&{\frac
{8045189}{871782912000}}{v}^{10}+{\frac {16009177}{10461394944000}}{v}
^{12}+O \left( {v}^{13} \right)
\end{eqnarray}
\begin{eqnarray}
b^{sd}_1=&&{\frac {5}{6}}-{\frac {1}{40}}{v}^{2}+{\frac {17}{2016}}{v}^{4}+{
\frac {1811}{1814400}}{v}^{6}+{\frac {13817}{79833600}}{v}^{8}+\nonumber \\&&{\frac
{12478951}{435891456000}}{v}^{10}+{\frac {24838031}{5230697472000}}{v}
^{12}+O \left( {v}^{13} \right)
\end{eqnarray}
\begin{eqnarray}
a(v)=&&-{\frac {1}{240}}{v}^{6}-{\frac {1}{2016}}{v}^{8}-{\frac {1}{11520}
}{v}^{10}-{\frac {2291}{159667200}}{v}^{12}-\nonumber \\&&{\frac {62879}{26417664000
}}{v}^{14}+O \left( {v}^{15} \right)
\end{eqnarray}

\begin{figure}[htbp]
\includegraphics[width=\textwidth]{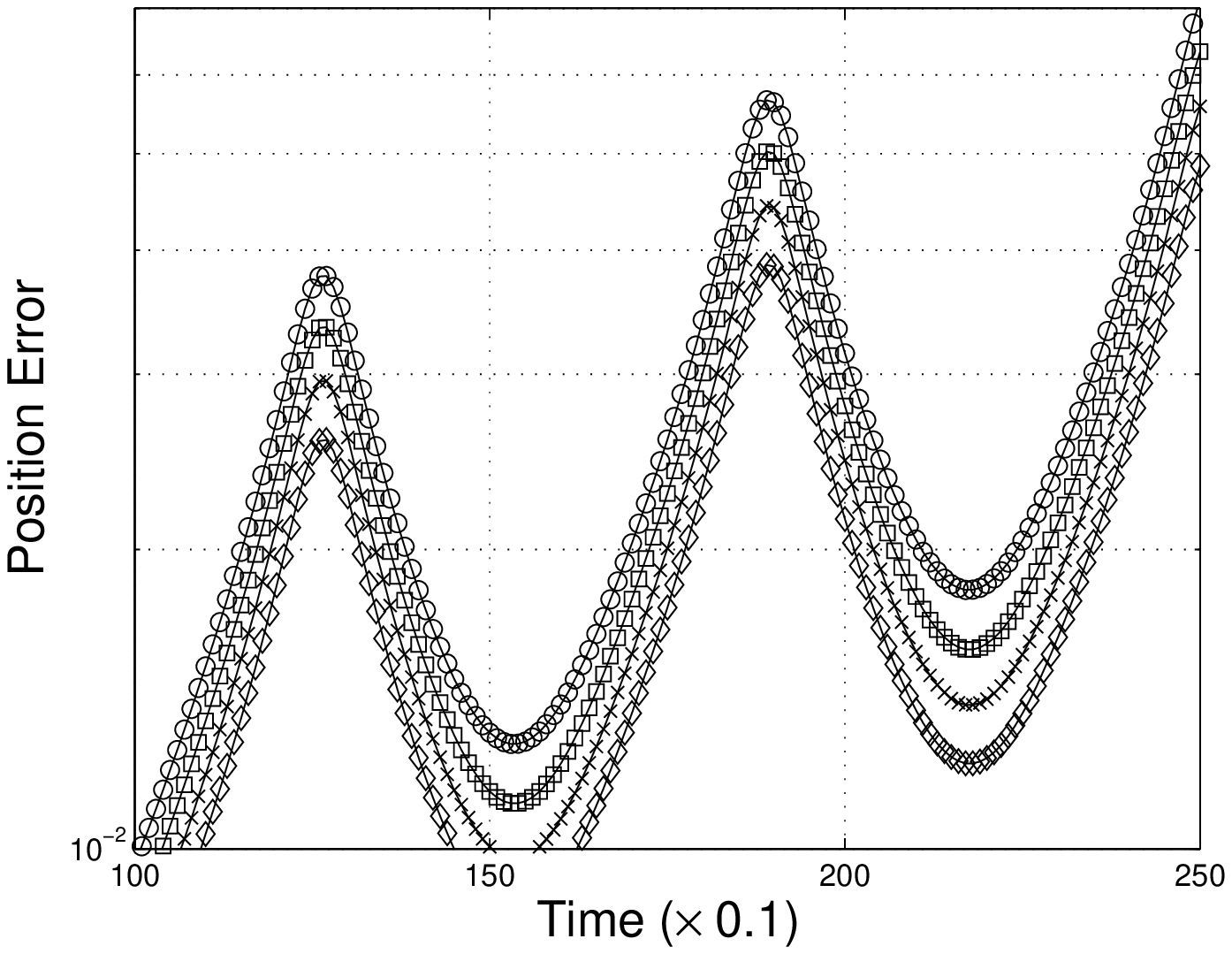}
\caption{Figure 1.}\label{fig_res_2b}
\end{figure}

\begin{figure}[htbp]
\includegraphics[width=\textwidth]{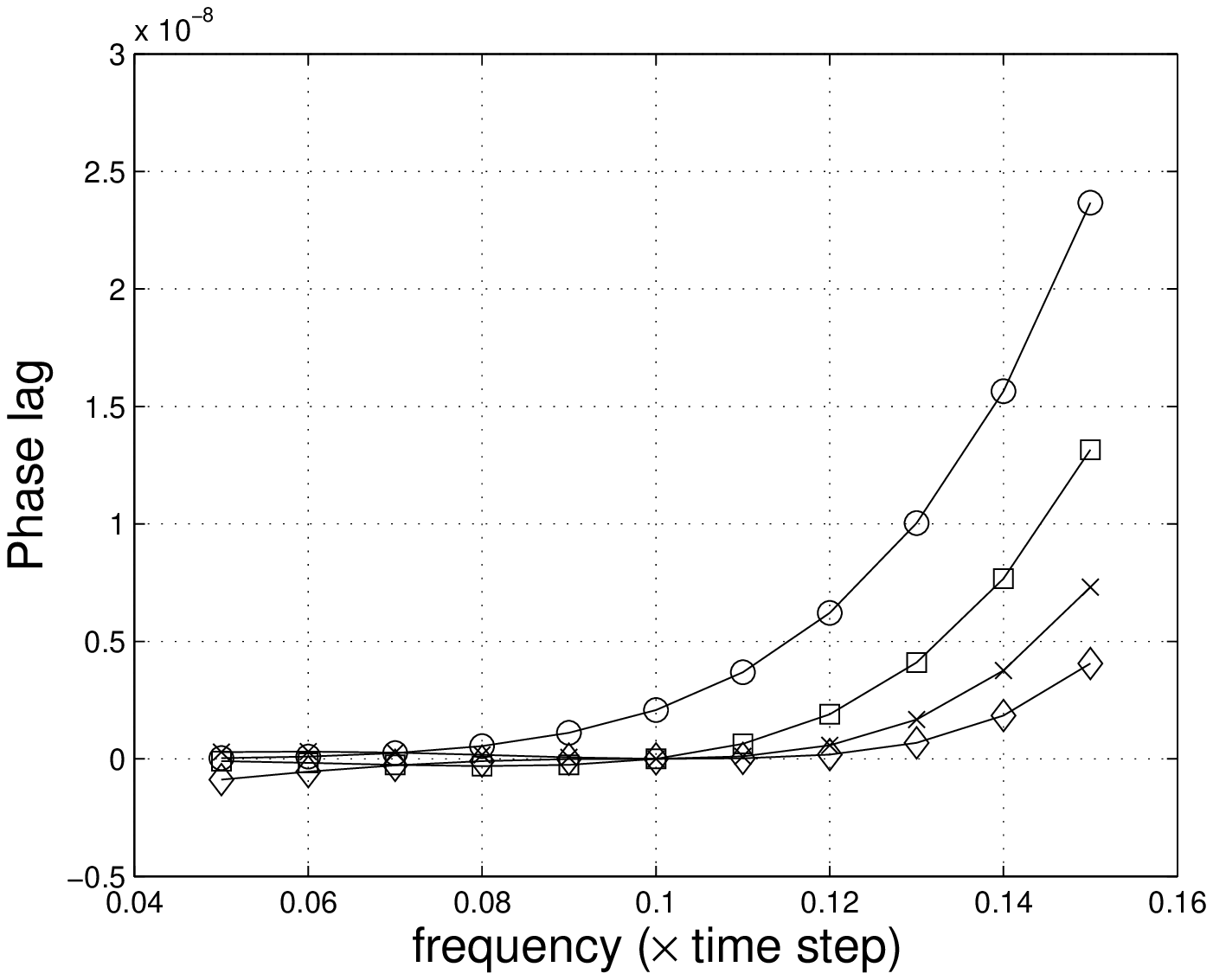}
\caption{Figure 2.}\label{fig_res_pl}
\end{figure}

The problem under test is the 2-body problem and the 2nd order equation is
\begin{equation}
\frac{d^2}{dt^2}\left( \begin{array}{l}y_1(t)\\y_2(t)\end{array}\right)=-\frac{1}{\sqrt{y_1^2(t)+y_2^2(t)}}\left(\begin{array}{l}y_1(t)\\y_2(t)\end{array}\right)
\end{equation}
with initial conditions
\begin{equation}
\left(\begin{array}{c}y_1(0) \\ y_2(0)\end{array}\right)=\left(\begin{array}{c}1-\epsilon \\ 0\end{array}\right),\;\left(\begin{array}{c}y_1(0) \\ y_2(0)\end{array}\right)'=\left(\begin{array}{c}0 \\ \sqrt{\frac{1+\epsilon}{1-\epsilon}}\end{array}\right)
\end{equation}
and exact solution
\begin{equation}
\left( \begin{array}{l}y_1(t)\\y_2(t)\end{array}\right)=\left( \begin{array}{l}cos(u)-\epsilon \\ \sqrt{1-\epsilon ^2}sin(u) \end{array} \right)
\end{equation}
and
\begin{equation}
u-\epsilon sin(u)-t=0
\end{equation}
In Figure \ref{fig_res_2b} the error in position calculation is plotted for the four methods. The working frequency in the trigonometric fitted methods has been estimated by $\omega=\frac{1}{(y_1^2+y_2^2)^\frac{3}{4}}$ \cite{anastassi_IJMPC_15_1_04}.
Finally in Figure \ref{fig_res_pl} the phase lag for the four methods is shown as a function of the frequency.

\section{Conclusions}
A new technique has been developed in order to improve the behavior of exponential (trigonometric) fitted numerical methods for the integration of oscillatory problems. The new technique is based on the vanishing of the first derivatives of the phase lag function, thus decreasing the sensitivity of the numerical method to frequency variations. Moreover, it has been shown that the new method becomes the classical one (the one who maximizes the algebraic order) when the working frequency tends to zero.

\bibliographystyle{plain}
\bibliography{../../../math}

\clearpage
\section*{List of Figures}
\begin{enumerate}
\item Mean error in position calculation in the two body problem for eccentricity $\epsilon=0.5$ and step size $h=0.1$. ($\circ$) is the algebraic fitted method, ($\boxdot$) is the trigonometric fitted method, ($\times$) is the trigonometric fitted method with the first derivative of the phase lag function equal to zero and ($\diamondsuit$) is the trigonometric fitted method with both first and second derivatives of the phase lag function equal to zero.
\item The phase lag of the four methods as a function of frequency. ($\circ$) is the algebraic fitted method, ($\boxdot$) is the trigonometric fitted method, ($\times$) is the trigonometric fitted method with the first derivative of the phase lag function equal to zero and ($\diamondsuit$) is the trigonometric fitted method with both first and second derivatives of the phase lag function equal to zero.
\end{enumerate}

\end{document}